\documentclass{gtart}

\def\ifplaintex{\expandafter\ifx\csname documentclass\endcsname\relax}

\def\gtp{{\mathsurround=0pt\it $\cal G\mskip-2mu$eometry \&\ 
$\cal T\!\!$opology $\cal P\!$ublications}}  % GT publications

\def\recd{{\small Received:\qua\receiveddate\ifx\reviseddate\relax
\else\qquad Revised:\qua\reviseddate\fi\par}} 

\def\lognumber#1{\def\thelognumber{#1}}
\def\volumenumber#1{\def\thevolumenumber{#1}}
\def\volumeyear#1{\def\thevolumeyear{#1}}
\def\papernumber#1{\def\thepapernumber{#1}}
\def\pagenumbers#1#2{\def\startpage{#1}\def\finishpage{#2}}
\def\published#1{\def\publishdate{#1}}

\def\received#1{\def\receiveddate{#1}}

\def\accepted#1{\def\accepteddate{#1}}

\def\asciiaddress#1{\def\theasciiaddress{#1}}
\def\asciiemail#1{\def\theasciiemail{#1}}

\long\def\asciiabstract#1{\long\def\theasciiabstract{#1}}

\let\\\par\let\thelognumber\relax\let\thevolumenumber\relax
\let\thepapernumber\relax\let\thevolumeyear\relax\let\startpage\relax
\let\finishpage\relax\let\publishdate\relax\let\receiveddate\relax
\let\reviseddate\relax\let\accepteddate\relax\let\theasciititle\relax
\let\theasciiauthors\relax\let\theasciiaddress\relax
\let\theasciiabstract\relax

\let\theasciiemail\relax

\ifplaintex
\font\logobig=cmssbx10 scaled 3836
\font\logomed=cmssbx10 scaled 2557
\else
\font\logobig=cmssbx10 scaled 4200
\font\logomed=cmssbx10 scaled 2800
\fi

\long\def\makeagttitle{   %%% start of definition of \makeagttitle
\count0=\startpage
\agt\hfill      %   Journal title (top left) 
\hbox to 45truept{\vbox to 0pt{\vglue -13truept{\logomed A\kern -.37em{\logobig 
T}\kern -.38em G}\vss}\hss}
\break
{\small Volume \thevolumenumber\ (\thevolumeyear)
\startpage--\finishpage\nl
Published: \publishdate}

\vglue .25truein

{\parskip=0pt\leftskip 0pt plus
1fil\def\\{\par\smallskip}{\Large\bf\thetitle}\par\medskip} \vglue
0.05truein

{\parskip=0pt\leftskip 0pt plus 1fil\def\\{\par}{\sc\theauthors}
\par\medskip}%
 
\vglue 0.03truein 

{\small\leftskip 25truept\rightskip 25truept{\bf Abstract}\stdspace\theabstract

{\bf AMS Classification}\stdspace\theprimaryclass
\ifx\thesecondaryclass\relax\else; \thesecondaryclass\fi\par
{\bf Keywords}\stdspace \thekeywords\par}\vglue 7truept

}   %%%% end of definition of \makeagttitle

\ifplaintex
\font\phead=cmsl9 scaled 950
\font\pnum=cmbx10 scaled 913
\font\pfoot=cmsl9 scaled 950
\headline{\vbox to 0pt{\vskip -4.5mm\line{\small\phead\ifnum
\count0=\startpage ISSN 1472-2739 (on-line) 1472-2747 (printed)
\hfill {\pnum\folio}\else\ifodd\count0\def\\{ }% 
\ifx\theshorttitle\relax\thetitle\else\theshorttitle\fi\hfill{\pnum\folio}
\else\def\\{ and }{\pnum\folio}\hfill\ifx\theshortauthors\relax\theauthors
\else\theshortauthors\fi\fi\fi}\vss}}
\footline{\vbox to 0pt{\vglue 0mm\line{\small\pfoot\ifnum\count0=\startpage
\copyright\ \gtp\hfill\else
\agt, Volume \thevolumenumber\ (\thevolumeyear)\hfill\fi}\vss}}
\else
\font=cmsl9 scaled 1050
\font\lnum=cmbx10 
\font=cmsl9 scaled 1050
\makeatletter
\def\@oddhead{{\small\lhead\ifnum\count0=\startpage ISSN 1472-2739 
(on-line) 1472-2747 (printed)\hfill {\lnum\number\count0}\else\ifodd\count0
\def\\{ }\ifx\theshorttitle\relax \thetitle \else\theshorttitle\fi\hfill
{\lnum\number\count0}\else\def\\{ and }{\lnum\number\count0}
\hfill\ifx\theshortauthors\relax 
\theauthors\else\theshortauthors\fi\fi\fi}}\def\@evenhead{\@oddhead}
\def\@oddfoot{\small\lfoot\ifnum\count0=\startpage\copyright\ \gtp\hfill\else
\agt, Volume \thevolumenumber\ (\thevolumeyear)\hfill\fi}
\def\@evenfoot{\@oddfoot}
\makeatother
\fi
\let\maketitlepage\makeagttitle
\let\makeshorttitle\maketitlepage
\let\maketitle\maketitlepage

\newwrite\gtoutfile
\long\gdef\makeheadfile{  %%% start of definition of \makeheadfile
{\def\\{, }\def\s{ }
\immediate\openout\gtoutfile head.xxx
\immediate\write\gtoutfile{To: math@arxiv.org}
\immediate\write\gtoutfile{Subject: put OR rep NNNNN:ppppp}
\immediate\write\gtoutfile{--text follows this line--}
\immediate\write\gtoutfile{Proxy-for: \ifx\theasciiauthors\relax
\theauthors\else\theasciiauthors\fi\s<\ifx\theasciiemail\relax\theemail\else\theasciiemail\fi>}
\immediate\write\gtoutfile{\noexpand\\}
\immediate\write\gtoutfile{Authors: \ifx\theasciiauthors\relax
\theauthors\else\theasciiauthors\fi}
{\def\\{ }\immediate\write\gtoutfile{Title: \ifx\theasciititle\relax
\thetitle\else\theasciititle\fi}}
\immediate\write\gtoutfile{Subj-class: GT or SG, GR etc}
\immediate\write\gtoutfile{MSC-class: \theprimaryclass\ifx\thesecondaryclass\relax\else, \thesecondaryclass\fi}
\immediate\write\gtoutfile{Journal-ref: Algebr. Geom. Topol. \thevolumenumber\s
(\thevolumeyear) \startpage-\finishpage}
\immediate\write\gtoutfile{Comments: Published by Algebraic and
Geometric Topology at}
\immediate\write\gtoutfile{\s\s\s  http://www.maths.warwick.ac.uk/agt/AGTVol\thevolumenumber/agt-\thevolumenumber-\thepapernumber.abs.html}
\immediate\write\gtoutfile{\noexpand\\}
\immediate\write\gtoutfile{}
\ifx\theasciiabstract\relax
\immediate\write\gtoutfile{\theabstract}\else
\immediate\write\gtoutfile{\theasciiabstract}\fi
\immediate\write\gtoutfile{}
\immediate\write\gtoutfile{\noexpand\\}
\immediate\write\gtoutfile{}
\immediate\closeout\gtoutfile}}  %%% end of definition of \makeheadfile

\def\maketitlepage{\makeagttitle\makeheadfile}
\let\makeshorttitle\maketitlepage
\let\maketitle\maketitlepage

\lognumber{16}
\volumenumber{3}
\volumeyear{2003}
\papernumber{16}
\published{15 June 2003}
\pagenumbers{511}{535}
\received{14 February 2003}
\accepted{19 May 2003}

\usepackage{amscd,amsmath,amssymb}
\usepackage{subfigure}
\usepackage[all,poly,knot]{xy}

\let\url\undefined
\usepackage{hyperref}

\newcommand{\arxiv}[1]
{\texttt{\href{http://front.math.ucdavis.edu/#1}{#1}}}
\newcommand{\MR}[1]
{\href{http://www.ams.org/mathscinet-getitem?mr=#1}{MR #1}}

\theoremstyle{plain}
\newtheorem{thm}[subsection]{Theorem}
\newtheorem{prop}[subsection]{Proposition}
\newtheorem{lem}[subsection]{Lemma}
\newtheorem{cor}[subsection]{Corollary}
\newtheorem{theo}{Theorem}
\newtheorem*{question}{Question}

\theoremstyle{definition}

\theoremstyle{remark}
\newtheorem{rem}[subsection]{Remark}
\newtheorem*{ack}{Acknowledgments}

\newcommand{\tsum}{\textstyle\sum\nolimits}
\newcommand{\abs}[1]{\left|#1\right|}
\newcommand{\set}[1]{\left\{#1\right\}}
\newcommand{\angl}[1]{\left<#1\right>}
\renewcommand{\b}[1]{\mathbf{#1}}

\newcommand{\A}{\mathcal{A}}
\newcommand{\B}{\mathcal{B}}

\newcommand{\C}{\mathbb{C}}
\newcommand{\F}{\mathbb{F}}
\newcommand{\K}{\mathbb{K}}
\newcommand{\Z}{\mathbb{Z}}
\newcommand{\II}{\mathbb{I}}
\newcommand{\CP}{{\mathbb{CP}}}

\newcommand{\T}{{\mathbf{T}}}
\newcommand{\FF}{{\mathbf{F}}}

\newcommand{\sA}{{\sf A}}
\newcommand{\sC}{{\sf C}}
\newcommand{\sM}{{\sf M}}

\DeclareMathOperator{\rank}{rank}
\DeclareMathOperator{\Aut}{Aut}
\DeclareMathOperator{\Hom}{Hom}
\DeclareMathOperator{\GL}{GL}
\DeclareMathOperator{\con}{cong}
\DeclareMathOperator{\ii}{i}
\DeclareMathOperator{\ord}{ord}
\DeclareMathOperator{\ab}{ab}
\DeclareMathOperator{\depth}{depth}
\DeclareMathOperator{\ch}{char}

\newcommand{\pmi}{\phantom{-}}
\newcommand{\dA}{{\mathsf{d}\mathcal{A}}}
\newcommand{\surj}{{\twoheadrightarrow}}
\newcommand{\cross}{{\times}}
\newcommand{\moo}{{U}}

\begin{document}

\title{Torsion in Milnor fiber homology}

\authors{Daniel C. Cohen\\Graham Denham\\Alexander I. Suciu}
\shortauthors{Cohen, Denham and Suciu}

\addresses{Department of Mathematics, Louisiana State University\\
Baton Rouge, LA 70803, USA\\\smallskip\\
Department of Mathematics, University of Western Ontario\\
London, ON  N6A 5B7, Canada\\\smallskip\\
Department of Mathematics, Northeastern University,
Boston, MA 02115, USA}
\asciiaddress{Department of Mathematics, Louisiana State University\\
Baton Rouge, LA 70803, USA\\
Department of Mathematics, University of Western Ontario\\
London, ON  N6A 5B7, Canada\\
Department of Mathematics, Northeastern University,
Boston, MA 02115, USA}

\email{\href{mailto:cohen@math.lsu.edu}{cohen@math.lsu.edu},
\href{mailto:gdenham@uwo.ca}{gdenham@uwo.ca},
\href{mailto:a.suciu@neu.edu}{a.suciu@neu.edu}}

\asciiemail{cohen@math.lsu.edu, gdenham@uwo.ca, a.suciu@neu.edu}

\urladdr{\href{http://www.math.lsu.edu/~cohen/}%
{http://www.math.lsu.edu/\char'176cohen},
\href{http://www.math.uwo.ca/~gdenham}%
{http://www.math.uwo.ca/\char'176gdenham},
\href{http://www.math.neu.edu/~suciu/}%
{http://www.math.neu.edu/\char'176suciu}}

\begin{abstract}
In a recent paper, Dimca and N\'emethi pose the problem of finding a
homogeneous polynomial $f$ such that the homology of the complement of
the hypersurface defined by $f$ is torsion-free, but the homology of
the Milnor fiber of $f$ has torsion.  We prove that this is indeed
possible, and show by construction that, for each prime $p$, there is 
a polynomial with $p$-torsion in the homology of the Milnor fiber.  
The techniques make use of properties of characteristic varieties of
hyperplane arrangements.
\end{abstract}

\asciiabstract{In a recent paper, Dimca and Nemethi pose the problem
of finding a homogeneous polynomial f such that the homology of the
complement of the hypersurface defined by f is torsion-free, but the
homology of the Milnor fiber of f has torsion.  We prove that this is
indeed possible, and show by construction that, for each prime p,
there is a polynomial with p-torsion in the homology of the Milnor
fiber.  The techniques make use of properties of characteristic
varieties of hyperplane arrangements.}

% 32S55 Milnor fibration; relations with knot theory
%       (Several complex variables and analytic spaces; Singularities)
% 14J70 Hypersurfaces
%       (Surfaces and higher-dimensional varieties)
% 32S22 Relations with arrangements of hyperplanes
%       (Several complex variables and analytic spaces; Singularities)
% 55N25 Homology with local coefficients, equivariant cohomology
%       (Algebraic topology; Homology and cohomology theories)


\maketitle

\section{Introduction}
\label{sec:intro}

Let $f\co (\C^{\ell+1},\b{0}) \to (\C,0)$ be a homogeneous polynomial. 
Denote  by $M=\C^{\ell+1} \setminus f^{-1}(0)$ the complement of the 
hypersurface defined by the vanishing of $f$, and let $F=f^{-1}(1)$ be
the Milnor fiber of the bundle map $f\co M\to \C^*$.  In
\cite[Question~3.10]{DN}, Dimca and N\'emethi ask the following.
\begin{question}
\label{quest:adan}
Suppose the integral homology of $M$ is torsion-free.  Is then
the integral homology of $F$ also torsion-free?
\end{question}

The Milnor fiber $F$ has the homotopy type of a finite,
$\ell$-dimensional CW-complex.  If $f$ has an isolated singularity at
$\b{0}$ (for example, if $\ell=1$), then $F$ is homotopic to a bouquet
of $\ell$-spheres, and so $H_*(F;\Z)$ is torsion-free.  The purpose of
this paper is to prove the following result, which provides a negative
answer to the above question, as soon as $\ell>1$.

\begin{theo}
\label{theo:mftors}
Let $p$ be a prime number, and let $\ell$ be an integer greater than
$1$.  Then there is a homogeneous polynomial $f_{p,\ell}\co 
\C^{\ell+1} \to \C$ for which $H_*(M;\Z)$ is torsion-free, but
$H_1(F;\Z)$ has $p$-torsion.
\end{theo}

Let $x_1,\dots,x_{\ell+1}$ be coordinates for $\C^{\ell+1}$.  The
theorem is proven by finding criteria for the construction of such
polynomials, then by explicitly exhibiting a family of $3$-variable
polynomials $f_p=f_p(x_1,x_2,x_3)$ with the desired properties, for
all primes $p$:
\begin{equation}
\label{eq:fp}
f_p = \begin{cases}
x_1x_2 (x_1^p-x_2^p)^{2} (x_1^p- x_3^p) (x_2^p- x_3^p),
&\text{if $p$ is odd,}
\\
x_1^2x_2(x_1^2-x_2^2)^3(x_1^2-x_3^2)^2(x_2^2-x_3^2),
&\text{if $p=2$.}
\end{cases}
\end{equation}
It then suffices to take
$f_{p,\ell}(x_1,\dots,x_{\ell+1})=f_p(x_1,x_2,x_3)$.

The above polynomials are all products of powers of linear factors,
and so define multi-arrangements of hyperplanes.  See \cite{OT} as a
general reference on arrangements.  For each prime $p$, the underlying
arrangement $\A_p$ is a deletion of the arrangement associated to the
complex reflection group $G(3,1,p)$, and has defining polynomial
$Q(\A_p)=x_1x_2(x_1^p-x_2^p)(x_1^p-x_3^p)(x_2^p-x_3^p)$.  As is well
known, for any hyperplane (multi)-arrangement, the homology groups of
the complement are finitely-generated and torsion-free.  Thus, Theorem
\ref{theo:mftors} is a consequence of the following result, which
identifies more precisely the torsion in the homology of the Milnor
fiber of the corresponding multi-arrangement.

\begin{theo}
\label{theo:h1fp}
Let $F_p=f_p^{-1}(1)$ be the Milnor fiber of the polynomial defined in
\eqref{eq:fp}.  Then:
\begin{equation*}
\label{eq:b1p}
H_1(F_p;\Z) = \begin{cases}
\Z^{3p+1}\oplus \Z_p \oplus T,
&\text{if $p$ is odd,}
\\
\Z^{3p+1} \oplus \Z_2\oplus \Z_2 ,
&\text{if $p=2$,}
\end{cases}
\end{equation*}
where $T$ is a finite abelian group satisfying $T \otimes \Z_q=0$ for
every prime $q$ such that $q\nmid 2(2p+1)$.
\end{theo}

The $p$-torsion in $H_1(F_p;\Z)$ is the smallest it can be (without
being trivial).  Indeed, if $H_*(M;\Z)$ is torsion-free, then an
application of the Wang sequence for the Milnor fibration $F\to M\to
\C^*$ shows that if the $2$-torsion summand of $H_1(F;\Z)$ is
non-trivial, then it must contain a repeated factor (compare
\cite[Prop.~3.11]{DN}).

The complement $M$ of a (central) arrangement of $n$ hyperplanes
admits a minimal cell decomposition, that is, a cell decomposition for
which the number of $k$-cells equals the $k$-th Betti number, for each
$k\ge 0$, see \cite{Ra02}, \cite{DP}.  On the other hand, it is not
known whether the Milnor fiber of a reduced defining polynomial for
the arrangement admits a minimal cell decomposition.  As noted in
\cite{Ra02}, this Milnor fiber does admit a cell decomposition with
$n\cdot b_k(U)$ cells of dimension $k$, where $U$ is the complement of
the projectivized arrangement.  Our results show that there exist
multi-arrangements for which the Milnor fiber $F$ admits no minimal
cell decomposition.  Indeed, by the Morse inequalities, the existence
of such a cell decomposition would rule out torsion in $H_*(F;\Z)$.

This paper is organized as follows.  Relevant results concerning
finite abelian covers, characteristic varieties, and Milnor fibrations
of multi-arrangements are reviewed in Sections \ref{sec:jump} and
\ref{sec:milnor}.  Criteria which insure that the homology of the
Milnor fiber of a multi-arrangement has torsion are established in
Section \ref{sec:construction}.  Multi-arrangements arising from
deletions of monomial arrangements are studied in Sections
\ref{sec:loci} and \ref{sec:MoreHomology}.  The proof of
Theorem \ref{theo:h1fp} is completed in Section \ref{sec:proofthm}.

\begin{ack}
Some of this work was carried out while the authors attended the March
2002 Mini-Workshop ``Cohomology Jumping Loci'' at the Mathematisches
Forschungs\-institut Oberwolfach.  We thank the Institute for its
hospitality, and for providing an exciting and productive mathematical
environment.

This research was supported by Louisiana Board of Regents grant
LEQSF(1999-2002)-RD-A-01 and by National Security Agency grant
MDA904-00-1-0038 (D.~Cohen), by a grant from NSERC of Canada 
(G.~Denham), and by NSF grant DMS-0105342 (A.~Suciu). 
\end{ack}

\section{Finite abelian covers and cohomology jumping loci}
\label{sec:jump}

We start by reviewing some basic facts about finite abelian covers,
and how to derive information about their homology from the
stratification of the character torus of the fundamental group by
cohomology jumping loci.  A more detailed treatment in the case of
line arrangements may be found in the survey \cite{Su01}.

\subsection{Homology of finite abelian covers}
\label{subsec:hfac}
Let $(X,x_0)$ be a based, connected space with the homotopy type of a
finite CW-complex, and let $G=\pi_1(X,x_0)$ be its fundamental group.
Let $Y$ be a finite, regular, abelian cover of $X$, with deck
transformation group $A$.  Finally, let $\K$ be a field, with
multiplicative group of units $\K^{\cross}$, and let
$\widehat{G}=\Hom(G,\K^{\cross})$ be the group of $\K$-valued
characters of $G$.

We shall assume that $\K$ is algebraically closed, and that
the characteristic of $\K$ does not divide the order of $A$.
With these assumptions, finitely-generated $\K[A]$-modules are
semisimple.  Since $A$ is abelian, irreducible representations are
one-dimensional, given by characters $\chi\co  A\to\K^{\cross}$.
By composing with the map $G\twoheadrightarrow A$, we obtain
one-dimensional $G$-modules denoted $\K_\chi$.

The lemma below is not new, and its proof can be found in various
special cases.  See \cite{lib93}, \cite{Sa95}, \cite{MShall} in the
context of $2$-complexes; \cite{cs1} in the context of cyclic covers
of complements of arrangements; and \cite{ber86} in an algebraic
setting.  For completeness, we will sketch a proof of the version
needed here.

\begin{lem}
\label{lem:Hdecomp}
Let $p\co  Y\to X$ be a finite, regular, abelian cover with
group of deck transformations $A$, and let $\K$ be an algebraically
closed field, with $\ch\K\nmid \abs{A}$.  Then
\begin{equation}
\label{eq:irredecomp}
H_*(Y;\K)\cong\bigoplus_{\chi\in \widehat{A}}H_*(X;\K_\chi),
\end{equation}
where $\K_\chi$ denotes the rank one local system given by lifting
a character $\chi\in\widehat{A}=\Hom(A,\K^{\cross})$
to a representation of $G=\pi_1(X,x_0)$.  Furthermore, the
direct summand indexed by a character $\chi$ is the corresponding
isotypic component of $H_*(Y;\K)$ as a $\K[A]$-module.
\end{lem}
\begin{proof}
The Leray spectral sequence of the cover $p\co  Y\to X$
degenerates to give an isomorphism
\begin{equation*}
H_*(Y;\K)\cong H_*(X;\K[A]),
\end{equation*}
where the action of $G$ on $\K[A]$ is induced by left-multiplication
of $G$ on $A=\pi_1(X)/p_*(\pi_1(Y))$.  That is, $H_*(Y;\K)$ is the
homology of $C_*(Y)\otimes_{\K[G]}\K[A]$, a chain complex of
$A$-modules under the right action of $A$.  By our assumptions on
$\K$, all $\K[A]$-modules are semisimple, so the group algebra of $A$
is isomorphic, as an $A$-module, to a direct sum of (one-dimensional)
irreducibles: $\K[A]\cong \bigoplus_{\chi\in\widehat{A}}\K_\chi$.
This decomposition into isotypic components commutes with
$\otimes_{\K[G]}$ and homology, yielding \eqref{eq:irredecomp}.
\end{proof}

\subsection{Characteristic varieties}
\label{subsec:cvs}
Assume that $H_1(X;\Z)=G^{\ab}$ is torsion-free and non-zero, and fix
an isomorphism $\alpha\co  G^{\ab} \to \Z^n$, where $n=b_1(X)$.  Let
$\K$ be an algebraically closed field.  The isomorphism $\alpha$
identifies the character variety $\widehat{G}=\Hom(G,\K^{\cross})$
with the algebraic torus $\b{T}(\K)=(\K^{\cross})^n$.

The cohomology jumping loci, or {\em characteristic varieties}, of $X$
are the subvarieties $\Sigma^q_d(X,\K)$ of the character torus defined
by
\begin{equation}
\label{eq:cjl}
\Sigma^q_d(X,\K) = \{\b{t}=(t_1,\dots,t_n) \in
(\K^{\cross})^n \mid \dim_\K H^q(X;\K_\b{t}) \ge d\},
\end{equation}
where $\K_{\b{t}}$ denotes the rank one local system given by the
composite $G \xrightarrow{\ab}G^{\ab} \xrightarrow{\alpha} \Z^n
\xrightarrow{\b{t}} \K^{\cross}$, and the last homomorphism sends
the $j$-th basis element to $t_j$.  For fixed $q>0$, these loci
determine a (finite) stratification
\[
(\K^{\cross})^n \supseteq\Sigma^q_1(X,\K)
\supseteq\Sigma^q_2(X,\K)\supseteq\cdots \supseteq \emptyset.
\]
Define the {\em depth} of a character $\b{t}\co  G\to\K^{\cross}$
relative to this stratification by
\begin{equation*}
\label{eq:depth}
\depth^q_{X,\K}(\b{t}) = \max \{d \mid \b{t} \in \Sigma^q_d(X;\K)\}.
\end{equation*}

The varieties $\Sigma^1_d(G,\K)$, the jumping loci for $1$-dimensional
cohomology of the Eilenberg-MacLane space $K(G,1)$, are particularly
accessible.  Indeed, these varieties are the determinantal varieties
of the Alexander matrix associated to a (finite) presentation of $G$,
see for instance \cite[Rem.~5.2]{MShall}.

Now assume that $H_2(X;\Z)$ is also torsion-free, and that the
Hurewicz homomorphism $h\co  \pi_2(X)\to H_2(X)$ is the zero map.
Then $H^2(X)=H^2(G)$, and this readily implies
$\Sigma^1_d(X,\K)=\Sigma^1_d(G,\K)$.  Thus, we may compute
$\depth_{\K}(\b{t}):=\depth^1_{X,\K}(\b{t})$ directly from the
Alexander matrix of~$G$.

\subsection{Finite cyclic covers}
\label{subsect:cyc}
Consider the case where $A=\Z_N$ is a cyclic group of order $N$.
Assume the characteristic of the field $\K$ does not divide $N$, so
that the homomorphism $\iota\co  \Z_N \to \K^{\cross}$ which sends a
generator of $\Z_N$ to a primitive $N$-th root of unity in $\K$ is an
injection.  For a homomorphism $\lambda\co  G\to \Z_N$, and an
integer $j>0$, define a character $\lambda^j\co  G\to \K^{\cross}$
by $\lambda^j(g)=\iota(\lambda(g))^j$.

Let $X$ be a finite CW-complex, with $H_1(X)$ and $H_2(X)$
torsion-free, and such that the Hurewicz map $h\co  \pi_2(X)\to
H_2(X)$ is trivial.  In view of the preceding discussion, Theorem~6.1
in \cite{MShall} applies as follows.

\begin{cor}
\label{cor:finind}
Let $p\co  Y\to X$ be a regular, $N$-fold cyclic cover, with
classifying map $\lambda\co  \pi_1(X)\surj \Z_N$.  Let $\K$ be an
algebraically closed field, with $\ch\K\nmid N$.  Then
\[
\dim_\K H_1(Y;\K) = b_1(X) + \sum_{1\neq k|N} \varphi(k)
\depth_\K\big(\lambda^{N/k}\big),
\]
where $\varphi$ is the Euler totient function.
\end{cor}

This result was first used in \cite{MShall} to detect $2$-torsion in
the homology of certain $3$-fold covers of the complement of the
deleted $\operatorname{B}_3$ arrangement (see \S \ref{subsect:peven}
below).  We will apply this result to Milnor fibrations in what
follows.

\section{Homology of the Milnor fiber of a multi-arrange\-ment}
\label{sec:milnor}

In this section, we review some facts concerning the Milnor fibration
of a complex (multi)-arrangement of hyperplanes, following \cite{cs1}
and \cite{de02}.

\subsection{Hyperplane arrangements}
\label{subsec:hyp}
Let $\A$ be a central arrangement of hyperplanes in $\C^{\ell+1}$.
The union of the hyperplanes in $\A$ is the zero locus of a polynomial
\begin{equation*}
\label{eq:defpoly}
f=Q(\A)=\prod_{H\in\A}\alpha_H \, ,
\end{equation*}
where each factor $\alpha_H$ is a linear form with kernel $H$.  Let
$\C^*\to \C^{\ell+1} \setminus\{\b{0}\} \to \CP^{\ell}$ be the Hopf
bundle, with fiber $\C^*=\C\setminus\{0\}$.  The projection map of
this (principal) bundle takes the complement of the arrangement,
$M=M(\A)=\C^{\ell+1}\setminus f^{-1}(0)$, to the complement $U$ of the
projectivization of $\A$ in $\CP^{\ell}$.  The bundle splits over $U$,
and so $M=U\times\C^*$.

It is well known that $U$ is homotopy equivalent to a finite
CW-complex (of dimension at most $\ell$), and that $H_*(U;\Z)$ is
torsion-free.  Furthermore, for each $k \ge 2$, the Hurewicz
homomorphism $h\co \pi_k(U)\to H_k(U)$ is the zero map, see
\cite{Ra97}.  Thus, the assumptions from \S\ref{subsect:cyc} hold for
$X=U$.

The fundamental group $\pi_1(M)$ is generated by meridian loops
(positively oriented linking circles) about the hyperplanes of $\A$.
The homology classes of these loops freely generate $H_1(M)=\Z^{n}$,
where $n=\deg(f)=\abs{\A}$.  We shall abuse notation and denote both a
meridian loop about hyperplane $H\in \A$, and its image in $\pi_1(U)$
by the same symbol, $\gamma_H$.  Note that these meridians may be
chosen so that $\prod_{H\in\A}\gamma_H$ is null-homotopic in $U$.  In
fact, $\pi_1(U)\cong\pi_1(M)/\langle \prod_{H\in\A}\gamma_H \rangle$,
and so $H_1(U)=\pi_1(U)^{\ab}=\Z^{n-1}$.

\subsection{The Milnor fibration}
\label{subsec:mfa}
As shown by Milnor, the restriction of $f\co  \C^{\ell+1}\to \C$ to
$M$ defines a smooth fibration $f\co  M\to\C^*$, with fiber
$F=f^{-1}(1)$ and monodromy $h\co  F\to F$ given by multiplication
by a primitive $n$-th root of unity in~$\C$.

The restriction of the Hopf map to the Milnor fiber gives rise to an
$n$-fold cyclic covering $F\to U$.  This covering is classified by the
epimorphism $\lambda\co  \pi_1(U) \surj \Z_n$ that sends all
meridians $\gamma_H$ to the same generator of $\Z_n$.  See \cite{cs1}
for details.

Now fix an ordering $\A=\set{H_1,H_2,\ldots,H_n}$ on the set of
hyperplanes.  Let $\b{a}=(a_1,a_2,\ldots,a_n)$ be an $n$-tuple of
positive integers with greatest common divisor equal to $1$.  We call
such an $n$-tuple a {\em choice of multiplicities for $\A$}.  The
(unreduced) polynomial
\begin{equation*}
\label{eq:multipoly}
f_{\b{a}}=\prod_{i=1}^n\alpha_{H_i}^{a_i}
\end{equation*}
defines a multi-arrangement $\A_{\b{a}}=\big\{H_1^{(1)},\dots,
H_1^{(a_1)},\ldots, H_n^{(1)},\dots, H_n^{(a_n)}\big\}$.  Note that
$\A_{\b{a}}$ has the same complement $M$, and projective complement
$U$, as $\A$, for any choice of multiplicities.  Let $f_{\b{a}}\co 
M\to \C^*$ be the corresponding Milnor fibration.  As we shall see,
the fiber $F_{\b{a}}=f_{\b{a}}^{-1}(1)$ does depend significantly on
$\b{a}$.

\subsection{Homology of the Milnor fiber}
\label{subsec:mfma}

Let $N=\sum_{i=1}^na_i$ be the degree of $f_{\b{a}}$, and let
$\Z_N=\langle g\mid g^N=1\rangle$ be the cyclic group of order $N$,
with fixed generator $g$.  As in the reduced case above, the
restriction of the Hopf map to $F_{\b{a}}$ gives rise to an $N$-fold
cyclic covering $F_{\b{a}} \to U$, classified by the homomorphism
$\lambda_{\b{a}}\co  \pi_1(U) \surj \Z_N$ which sends the meridian
$\gamma_i$ about $H_i$ to $g^{a_i}$.

For any field $\K$, let $\tau\co  (\K^\cross)^n \to \K^\cross$ be
the map which sends an $n$-tuple of elements to their product.  Since
the meridians $\gamma_i$ may be chosen so that $\prod_{i=1}^n
\gamma_i=1$, if $\b{s} \in (\K^\cross)^n$ satisfies $\tau(\b{s})=1$,
then $\b{s}$ gives rise to a rank one local system on $U$, compare
\S\ref{subsec:cvs}.  We abuse notation and denote this local system by
$\K_{\b{s}}$.

Suppose that $\K$ is algebraically closed, and $\ch\K$ does not divide
$N$.  Then there is a primitive $N$-th root of unity $\xi\in \K$.  Let
$\b{t}\in (\K^{\cross})^n$ be the character with $t_i=\xi^{a_i}$, for
$1\leq i\leq n$.  Note that $\tau(\b{t})=1$.  Let $h_{\b{a}}\co 
F_{\b{a}}\to F_{\b{a}}$ be the geometric monodromy of the Milnor
fibration $f_{\b{a}}\co  M \to \C^*$, given by multiplying
coordinates in $\C^{\ell+1}$ by a primitive $N$-th root of unity in
$\C$.  The action of the algebraic monodromy $(h_{\b{a}})_*\co 
H_*(F_{\b{a}};\K) \to H_*(F_{\b{a}};\K)$ coincides with that of the
deck transformations of the covering $F_{\b{a}} \to U$.
Lemma~\ref{lem:Hdecomp} yields the following.

\begin{lem}
\label{lem:Fsplits}
With notation as above, we have
\begin{equation*}
H_*(F_{\b{a}};\K)=\bigoplus_{k=0}^{N-1}H_*(U;\K_{\b{t}^k}).
\end{equation*}
Furthermore, the $k$-th summand is an eigenspace for $(h_{\b{a}})_*$
with eigenvalue $\xi^k$.
\end{lem}

The next lemma appeared in \cite{de02} in the complex case.  For
convenience, we reproduce the proof in general.

\begin{lem}
\label{lem:mults}
Let $\K$ be an algebraically closed field, and let $\b{s}\in
(\K^{\cross})^n$ be an element of finite order, with $\tau(\b{s})=1$.
Then there exists a choice of multiplicities $\b{a}$ for $\A$ so that
$H_q(U;\K_\b{s})$ is a monodromy eigenspace of $H_q(F_\b{a};\K)$.
\end{lem}

\begin{proof}
Let $\zeta\in\K$ be a primitive $k$-th root of unity, where $k$ is the
order of $\b{s}$.  Then, for each $2\leq i\leq n$, there is an integer
$1\leq a_i\leq k$ such that $s_i=\zeta^{a_i}$.  By choosing either
$1\leq a_1\leq k$ or $k+1\leq a_1 \leq 2k$ suitably, we can arrange
that the sum $N=\sum_{i=1}^n a_i$ is not divisible by $p=\ch\K$, if
$p>0$.  Since $\b{s}$ and $\zeta$ both have order $k$, we have
$\gcd\set{a_1,\dots ,a_n}=1$.  Since the product of the coordinates of
$\b{s}$ is $1$, the integer $k$ divides $N$.

By insuring $p\nmid N$, there is an element $\xi\in\K$ for which
$\xi^{N/k}=\zeta$.  By construction, $\b{a}=(a_1,\ldots,a_n)$ is a
choice of multiplicities for which $\b{s}=\b{t}^{N/k}$ in the
decomposition of Lemma~\ref{lem:Fsplits}, so $H_q(U;\K_\b{s})$ is a
direct summand of $H_q(F_\b{a};\K)$.
\end{proof}

\begin{rem}
The choice of multiplicities $\b{a}$ in Lemma~\ref{lem:mults} is not
unique.  As above, write $s_i=\zeta^{a_i}$ for integers $a_i$, where
$\zeta$ is a $k$-th root of unity and $1\leq a_i\leq k$.  Let
$\b{a}=(a_1,\ldots,a_n)$.  Then $H_q(U;\K_\b{s})$ is also a monodromy
eigenspace of $F_\b{b}$ if $\b{b}=\b{a}+\lambda$, for all
$\lambda\in(k\Z)^n$ for which satisfy $b_i>0$ for each $i$ and, if
$p>0$, $p\nmid\sum_{i=1}^n b_i$.
\end{rem}

\section{Translated tori and torsion in homology}
\label{sec:construction}

\subsection{Characteristic varieties of arrangements}
\label{subsec:cvarr}
Let $\A=\set{H_1,\ldots,H_n}$ be a central arrangement in
$\C^{\ell+1}$.  Let $M$ denote its complement, and $U$ the complement
of its projectivization.  Then the restriction of the Hopf fibration
$\C^*\rightarrow M\rightarrow U$ induces an isomorphism
$\pi_1(U)\cong\pi_1(M)/\angl{\prod_{i=1}^n\gamma_i}$, as in the
previous section.  For this reason, although the rank of
$\pi_1(U)^{\rm ab}$ is $n-1$, we shall regard the characteristic
varieties of $U$ as embedded in the character torus of $\pi_1(M)$:
\begin{equation}
\label{eq:sigmau}
\Sigma^q_d(U,\K)=\set{\b{t}\in\ker\tau\cong (\K^\times)^{n-1} \mid
\dim_\K H^q(U;\K_\b{t}) \ge d},
\end{equation}
(compare with \eqref{eq:cjl}), where, as above,
$\tau\co (\K^\times)^n\rightarrow\K^\times$
is the homomorphism given by $\tau(t_1,\ldots,t_n)=t_1\cdots t_n$.

\begin{prop}  \label{prop:cvMandU}
For $q\geq1$ and $d\geq0$,
\begin{equation*}
\Sigma^q_d(M,\K)=\bigcup_{j=0}^d \Sigma^q_{d-j}(U,\K)\cap
\Sigma^{q-1}_j(U,\K).
\end{equation*}
In particular, for $\b{t}\in \ker\tau$, we have
\begin{equation*}
\depth^q_{M,\K}(\b{t})=
\depth^q_{U,\K}(\b{t})+\depth^{q-1}_{U,\K}(\b{t}).
\end{equation*}
\end{prop}
\begin{proof}
Let $\K_\b{t}$ be the local system on $M$ corresponding to
$\b{t}\in(\K^\times)^n$.  There is an induced local system
$i^*\K_\b{t}$ on $\C^*$, with monodromy $\tau(\b{t})$, where
$i\co \C^* \to M$ is the inclusion of the fiber in the
Hopf bundle $\C^* \to M \to U$.  Fix a section $s\co  U \to M$
of this trivial bundle, and let $s^*\K_\b{t}$ be the induced
local system on $U$.  Recall that we denote this local system
by $\K_\b{t}$ in the case where $\tau(\b{t})=1$.
To prove the Proposition, it suffices to show that,
for each $q\ge 1$,
\[
H^q(M;\K_\b{t})=
\begin{cases}
0, &\text{if $\tau(\b{t}) \neq 1$}, \\
H^q(U;\K_\b{t}) \oplus H^{q-1}(U;\K_\b{t}),
&\text{if $\tau(\b{t})=1$}.
\end{cases}
\]

Let $C_\bullet(\widetilde{M})$ and $C_\bullet(\widetilde{U})$
be the chain complexes of the universal covers of $M$ and $U$,
viewed as modules over the group rings of $\bar{G}=\pi_1(M)$ and
$G=\pi_1(U)$, respectively.  Then the cohomology of $M$ with
coefficients in $\K_\b{t}$ is (by definition) the cohomology
of the complex
$\bar{\sC}^\bullet=\Hom_{\Z{\bar{G}}}(C_\bullet(\widetilde{M}),\K)$,
where the $\Z{\bar{G}}$-module structure on $\K$ is given by the
representation $\bar{G} \xrightarrow{\ab} \bar{G}^{\ab}
\xrightarrow{\alpha} \Z^n \xrightarrow{\b{t}} \K^{\cross}$.
Similarly, $H^*(U;s^*\K_\b{t})$ is the cohomology of the complex
$\sC^\bullet=\Hom_{\Z{G}}(C_\bullet(\widetilde{U}),\K)$.
Denote the boundary maps of the complexes $\bar{\sC}^\bullet$ and
$\sC^\bullet$ by $\Delta^\bullet$ and $\delta^\bullet$, respectively.

Multiplication by $1-\tau(\b{t})$ gives rise to a chain map
$\sC^\bullet \to \sC^\bullet$.  Since $M=U \times \C^*$ is
a product, and the monodromy of the induced local system
$i^*\K_\b{t}$ on $\C^*$ is $\tau(\b{t})$, the complex
$\bar{\sC}^\bullet$ may be realized as the mapping cone
of this chain map.  Explicitly, we have
$\bar{\sC}^q = \sC^q \oplus \sC^{q-1}$, and $\Delta^q\co 
\sC^q \oplus \sC^{q-1} \to \sC^{q+1} \oplus \sC^{q}$ is given by
\[
\Delta^q(x,y)=\bigl(\delta^q(x),\delta^{q-1}(y)+
(-1)^q(1-\tau(\b{t}))\cdot x\bigr).
\]

If $\tau(\b{t}) \neq 1$, it is readily checked that the complex
$\bar{\sC}^\bullet$ is acyclic.  If $\tau(\b{t})=1$, it follows
immediately from the above description of the boundary map
$\Delta^\bullet$ that
$H^q(\bar{\sC}^\bullet) \cong H^q(\sC^\bullet) \oplus
H^{q-1}(\sC^\bullet)$ for each $q$.
\end{proof}

Now let $\dA$ be the decone of $\A$ with respect to one of the
hyperplanes (which, after a linear change of variables, may be assumed
to be a coordinate hyperplane).  The complement, $M(\dA)$, in
$\C^\ell$ is diffeomorphic to the complement $U$ of the
projectivization of $\A$.  An isomorphism $\pi_1(U)\rightarrow
\pi_1(M(\dA))$ is obtained by deleting the meridian corresponding to
the deconing hyperplane.  Let $\pi\co  (\K^{\cross})^n \to
(\K^{\cross})^{n-1}$ be the map that forgets the corresponding
coordinate.  Then $\pi$ induces a bijection
$\pi_{\sharp}\co \Sigma^q_d(U,\K)\rightarrow \Sigma^q_d(M(\dA),\K)$.

If $\b{s}$ is a nontrivial character, then $H^0(U,\K_\b{s})=0$ and
$\depth^1_{U,\K}(\b{s})<n-1$.  Consequently, as shown in \cite{CScv}
using properties of Fitting ideals, for $q=1$ and $d<n$, the above
proposition simplifies to:
\begin{equation}
\label{eq:decone}
\Sigma^1_d(M(\A),\K) = \set{ \b{t} \in (\K^{\cross})^n \mid
\pi(\b{t}) \in \Sigma^1_d(M(\dA),\K) \text{ and } \tau(\b{t})=1}.
\end{equation}

Each irreducible component of $\Sigma^q_d(\moo,\C)$ (resp.,
$\Sigma^q_d(M(\dA),\C)$) is a torsion-translated subtorus of the
algebraic torus $\T(\C)=(\C^{\cross})^n$, see \cite{ara97}.  That is,
each component of $\Sigma_d^q(\moo,\C)$ is of the form $gT$, where $T$
is a subgroup of $\T(\C)$ isomorphic to a product of $0$ or more
copies of $\C^{\cross}$, and $g\in\T(\C)$ is of finite order.  Recall
that every algebraic subgroup of $\T(\K)$ can be written as the
product of a finite group with a subtorus~\cite[p.~187]{norbook}.  If
the order of an element $g\in\T(\K)$ is finite, we will denote its
order by $\ord(g)$.

\subsection{Jumping loci and the Milnor fibration}
\label{subsec:jumpmf}
Write $H_i=\ker(\alpha_i)$ and let
$f_{\b{a}}=\prod_{i=1}^n\alpha_i^{a_i}$ be the polynomial of degree
$N=\sum_{i=1}^n a_i$ corresponding to a choice of multiplicities
$\b{a}=(a_1,\dots,a_n)$ for $\A$.  Recall that $F_\b{a}$, the Milnor
fiber of $f_{\b{a}}\co  M\to \C^*$, is the regular, $N$-fold cyclic
cover of $U$ classified by the homomorphism $\lambda_{\b{a}}\co 
\pi_1(U)\surj \Z_N$ given by $\lambda_{\b{a}}(\gamma_i)=g^{a_i}$.
Recall also that $b_1(U)=\abs{\A}-1=n-1$.  From
Corollary~\ref{cor:finind}, we obtain the following.

\begin{thm}
\label{thm:hommf}
Let $\K$ be an algebraically closed field, with $\ch\K\nmid N$.  Then
\[
\dim_\K H_1(F_\b{a};\K) = n-1 + \sum_{1\neq k|N} \varphi(k)
\depth_\K\big(\lambda_{\b{a}}^{N/k}\big).
\]
\end{thm}

\subsection{Jumping loci in different characteristics}
\label{subsec:jumpchar}

Our goal for the rest of this section is to show that if a translated
torus $gT$ is a positive-dimensional component of a characteristic
variety $\Sigma_d^q(\moo,\C)$, but $T$ itself is not a component, then
there exist choices of multiplicities $\b{a}$ for which
$H_q(F_\b{a};\Z)$ has integer torsion (Theorem~\ref{th:transl}).  In
fact, we will describe how to choose such exponents explicitly, and
give a more general criterion for the existence of torsion
(Theorem~\ref{th:jumps}).

We start by comparing representations of the fundamental group over
fields of positive characteristic with those over $\C$.  Let $\zeta$
be a root of unity, and denote by $\Z[\zeta]$ the ring of cyclotomic
integers.

\begin{lem}
\label{lem:c}
Let $i\co \Z[\zeta]\rightarrow\C$ and
$j\co \Z[\zeta]\rightarrow\K$ be ring homomorphisms, and assume that
$i$ is an injection.  For any $\b{t}\in(\Z[\zeta]^{\cross})^n$ with
$\tau(\b{t})=1$, let $i_*\b{t}$ and $j_*\b{t}$ denote the images of
$\b{t}$ in $\T(\C)$ and $\T(\K)$, respectively.  Then
\begin{equation*}
\dim_\C H_q(U;\C_{i_*\b{t}})\leq
\dim_\K H_q(U;\K_{j_*\b{t}}).
\end{equation*}
\end{lem}
\begin{proof}
Since the character $\b{t}$ satisfies $\tau(\b{t})=1$, it gives rise
to a homomorphism $\psi\co  \Z{G}\rightarrow\Z[\zeta]$, where
$G=\pi_1(U)$ and $\Z{G}$ is the integral group ring.  Let
$K_*=C_*(\widetilde{U})\otimes_{\psi}\Z[\zeta]$ denote the
corresponding tensor product of the chain complex of the universal
cover of $U$ with $\Z[\zeta]$, a chain complex of
$\Z[\zeta]$-modules.
Then the homology groups under comparison are just those of
$K_*\otimes_{i_*\b{t}}\C$ and $K_*\otimes_{j_*\b{t}}\K$,
respectively.
Since the first map $i$ is flat, the inequality follows.
\end{proof}

\begin{lem}\label{lem:criteria}
Given an arrangement $\A$ and positive integers $q$, $d$,
the following two statements are equivalent.
\begin{enumerate}
\item
The characteristic variety $\Sigma_d^q(\moo,\C)$ contains an element
$g$ of finite order for which the cyclic subgroup
$\angl{g}\not\subseteq\Sigma_d^q(\moo,\C)$.  Moreover, there exists
$h\in\angl{g}\setminus\Sigma_d^q(\moo,\C)$ and a prime $p$ with
$p\mid\ord(g)$ but $p\nmid\ord(h)$.
\item
There exist $\b{s},\b{t}\in\T(\C)$, a prime $p$, and integer
$r\geq 1$ for which
\begin{enumerate}
\item[\rm(a)] $\depth^q_{\moo,\C}(\b{t}) < \depth^q_{\moo,\C}(\b{s})=d$;
\item[\rm(b)] $\ord(\b{s}\b{t}^{-1})=p^r$;
\item[\rm(c)] $p\nmid\ord(\b{t})$.
\end{enumerate}
\end{enumerate}
\end{lem}
\begin{proof}
$(1)\Rightarrow(2)$: Write $\angl{g}\cong\bigoplus_{i=1}^m
\Z/(p_i^{r_i}\Z)$, where the primes $p_1,p_2,\ldots,p_m$ are all
distinct.  For each $h\in\angl{g}$, define an $m$-tuple $\nu(h)$ as
follows: for $1\leq i\leq m$, let $\nu(h)_i=a_i$, where the projection
of $h$ to $\Z/(p_i^{r_i}\Z)$ has order $p_i^{a_i}$.  Clearly $0\leq
a_i\leq r_i$.

Let $S$ consist of those elements $h\in\angl{g}$ for which
$h\in\Sigma^q_d(\moo,\C)$.  Since characteristic varieties are closed
under cyclotomic Galois actions, two elements $h_1,h_2\in\angl{g}$ of
the same order are either both in $S$ or both not in $S$.  By
reordering the $p_i$'s, our hypothesis states that there exists
$h\notin S$ with $\nu(h)=(a_1,\ldots,a_j,0,\ldots,0)$, for some
nonzero integers $a_1,a_2,\ldots,a_j$, where $j<m$.  Choose $h \notin
S$ of this form for which $j$ is minimal.  Since $\b1\in S$ and
$\nu(\b1)=(0,0,\ldots,0)$, we have $j\ge 1$.
Then for some $h'\in\angl{g}$ of order $p_j^{r_j-a_j}$, we have
$\nu(hh')=(a_1,\ldots,a_{j-1},0,0,\ldots,0)$.  By minimality, $hh'\in
S$.  Then the pair of $\b{t}=h$ and $\b{s}=hh'$ together with $p=p_r$,
$r=r_j-a_j$ satisfy the conditions (2).

$(2)\Rightarrow(1)$: Let $g=\b{s}$, $h=\b{t}$, and $h'=gh^{-1}$.  By
hypothesis, $\ord(hh')=\ord(h)\ord(h')$, from which it follows that
$\angl{g}=\angl{hh'}=\angl{h,h'}$.  In particular, $h\in\angl{g}$, but
by (a), $h\notin\Sigma^q_d(\moo,\C)$.
\end{proof}

\subsection{Torsion jumps} \label{subsec:tj}
Once again, let $\K$ be an algebraically closed field of positive
characteristic $p$.

\begin{thm}\label{th:jumps}
If $\A$ is an arrangement for which the characteristic variety
$\Sigma^q_d(\moo,\C)$ satisfies one of  the equivalent conditions
of Lemma~\ref{lem:criteria}, then 
\[
\dim_\K H_q(U;\K_{\b{t}})\geq d.
\]
\end{thm}
\begin{proof}
Let $k=\ord(\b{t})$; from condition $(2)$, parts (b) and (c), we have
$\ord(\b{s})=p^rk$.  Let $\zeta$ be a root of unity in $\C$ of order
$p^rk$, so that $\b{s},\b{t}\in(\Z[\zeta]^{\cross})^n$.  Let
$j\co \Z[\zeta]\rightarrow\K$ be given by choosing a $k$-th root of
unity $j(\zeta)$ in $\K$.  Since $\ord(\b{s}\b{t}^{-1})$ is a power of
the characteristic of $\K$, we have $j_*g=j_*h$.  Then
\[
\dim_\K H_q(U;\K_{j_*\b{t}})=
\dim_\K H_q(U;\K_{j_*\b{s}}) \ge d  \, ,
\]
by condition (2)(a) and Lemma~\ref{lem:c}.
\end{proof}

\begin{cor}
\label{cor:tor}
Suppose $\A$ is an arrangement for which the characteristic variety
$\Sigma^q_d(\moo,\C)$ satisfies the equivalent conditions of
Lemma~\ref{lem:criteria}.  Then there is a choice of multiplicities
$\b{a}$ for $\A$ for which the group $H_q(F_\b{a};\Z)$ contains
$p$-torsion elements.
\end{cor}
\begin{proof}
Assume that $\b{t}\in\T(\C)$ satisfies condition (2)(a) of
Lemma~\ref{lem:criteria}.  Then, since $\b{t}\not\in
\Sigma_d^q(\moo,\C)$, we have $\dim_\C H_q(U;\C_{\b{t}})< \dim_\K
H_q(U;\K_{\b{t}})$.  Lemma~\ref{lem:mults} implies that there exists a
choice of multiplicities $\b{a}$ for which $H_q(U;\K_{\b{t}})$ and
$H_q(U;\C_{\b{t}})$ are monodromy eigenspaces.  Using
Lemmas~\ref{lem:Fsplits} and \ref{lem:c}, with one of the inequalities
being strict, we find that $\dim_\C H_q(F_\b{a};\C)<\dim_\K
H_q(F_\b{a};\K)$.  The result follows.
\end{proof}

The following statement is a special case of Theorem~\ref{th:jumps}
that applies to some specific behavior observed in characteristic
varieties (see \cite{suciu02} and \cite{triples}).  In particular, we
will use it in what follows to find torsion for our family of
examples.

\begin{thm}
\label{th:transl}
Let $\b{s}T$ be a component of $\Sigma^q_d(\moo,\C)$, where $T$ is a
subtorus of $\T(\C)$ and $\b{s}$ is a finite-order element in
$\T(\C)$.  Suppose that $T\not\subseteq\Sigma^q_d(\moo,\C)$.  Then
there exist choices of multiplicities $\b{a}$ for $\A$ for which the
group $H_q(F_\b{a};\Z)$ has $p$-torsion, for some prime $p$ dividing
$\ord(\b{s})$.
\end{thm}

\begin{proof}
First, note that $T$ is positive-dimensional, since $\b{1}$ is
contained in all non-empty characteristic varieties.  Since $T$ is not
contained in $\Sigma^q_d(\moo,\C)$, there exist infinitely many
finite-order elements $h\in T$ for which
$h\not\in\Sigma^q_d(\moo,\C)$.  (In fact, for each sufficiently large
integer $k$, there exist elements $h$ with $\ord(h)=k$ and
$h\notin\Sigma^q_d(\moo,\C)$.)

Choose any element $h$ as above, of order relatively prime to that of
$\b{s}$, and let $u=h^r$ for an $r$ for which $u^{\ord(\b{s})}=h$.
Let $g=\b{s}u$.  Then, by construction, $g$ and $h$ satisfy the first
condition of Lemma~\ref{lem:criteria}.  By Corollary~\ref{cor:tor},
$H_q(F_\b{a};\Z)$ has torsion of order $p$ for those $\b{a}$ given by
Lemma~\ref{lem:mults}.
\end{proof}

\section{Deletions of monomial arrangements}
\label{sec:loci}

Now we turn to a detailed study of arrangements obtained by deleting a
hyperplane from a monomial arrangement.  Using results from
\cite{suciu02} and \cite{triples}, we check that these arrangements
satisfy the hypotheses of Theorem \ref{th:transl}.  Hence, there are
corresponding multi-arrangements whose Milnor fibers have torsion in
homology.

\subsection{Fundamental group of the complement}
\label{subsect:pi1mono}
Let $\A_p$ be the arrangement in $\C^3$ defined by the homogeneous
polynomial
$Q(\A_p)=x_1^{}x_2^{}(x_1^p-x_2^p)(x_1^p-x_3^p)(x_2^p-x_3^p)$.  This
arrangement is obtained by deleting the hyperplane $x_3=0$ from the
complex reflection arrangement associated to the full monomial group
$G(3,1,p)$.

The projection $\C^3\to\C^2$ defined by $(x_1,x_2,x_3) \mapsto
(x_1,x_2)$
restricts to a bundle map $M(\A_p) \to M(\B)$, where $\B$ is defined
by $Q(\B)=x_1x_2(x_1^p-x_2^p)$.  The fiber of this bundle is the
complex line with $2p$ points removed.  Thus, $\A_p$ is a fiber-type
arrangement, with exponents $(1,p+1,2p)$.  Hence, the fundamental
group
$G(\A_p)=\pi_1(M(\A_p))$ may be realized as a semidirect product
\begin{equation}
\label{eq:extension}
G(\A_p) = \FF_{2p} \rtimes_\alpha G(\B),
\end{equation}
where $\FF_{2p}=\pi_1(\C\setminus \{2p \text{ points}\})$ is free on
$2p$ generators corresponding to the hyperplanes defined by
$(x_1^p-x_3^p)(x_2^p-x_3^p)$, and $G(\B)\cong \FF_{p+1} \times \Z$ is
the fundamental group of $M(\B)$.

The monodromy $\alpha\co  G(\B) \to \Aut(\F_{2p})$ which defines the
semidirect product structure \eqref{eq:extension} factors as
$G(\B) \xrightarrow{\eta} P_{2p} \hookrightarrow \Aut(\FF_{2p})$,
where the inclusion of the pure braid group on $2p$ strands $P_{2p}$
in $\Aut(\FF_{2p})$ is given by the restriction of the Artin
representation.  The ``braid monodromy'' $\eta\co  G(\B) \to P_{2p}$
may be determined using the techniques of \cite{CSbm}, \cite{CScc},
and \cite{mono}.  In fact, this map may be obtained by an appropriate
modification of the calculation in \cite[\S2.2]{mono} of the braid
monodromy of the full monomial arrangement defined by $x_3 Q(\A_p)$,
which we now carry out.

\subsection{Braid monodromy}
\label{subsect:bmono}
Fix a primitive $p$-th root of unity $\xi \in \C$.  Let $B_{2p}$ be
the full braid group on $2p$ strands, and let $\sigma_i$, 
$1\le i \le 2p-1$, be the standard generators.  The indices of the
strands correspond to the hyperplanes $H_{i3:r}=\ker(x_i-\xi^r x_3)$
and the generators $y_1,\dots,y_{2p}$ of $\FF_{2p}$, as indicated
below:
\[
\begin{matrix}
\hfill \text{strand \#}
& 1 & 2 & \cdots  & p & p+1 & p+2 & \cdots  & 2p \\
\text{hyperplane} & H_{13:p} & H_{13:p-1} & \cdots
& H_{13:1} &   H_{23:p} & H_{23:p-1} & \cdots & H_{23:1}\\
\hfill\text{generator} & y_1 & y_2 & \cdots & y_p
& y_{p+1} & y_{p+2} & \cdots & y_{2p}
\end{matrix}
\]

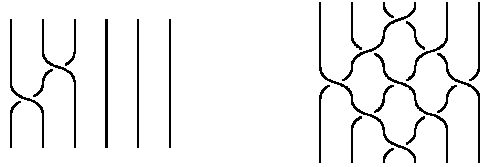
\begin{figure}%[ht]
\setlength{\unitlength}{0.7cm}
\subfigure{%
\label{fig:braid1}%
\begin{minipage}[t]{0.35\textwidth}
\begin{picture}(4,2.5)(-4.8,-0.7)
\xy
0;/r1pc/:
,{\vtwist~{(0,3)}{(0,1)}{(0,3)}{(0,1)}}
,{\vtwist~{(1,3)}{(1,2)}{(1,3)}{(1,2)}}
,{\vtwist~{(2,3)}{(2,2)}{(2,3)}{(2,2)}}
,{\htwist~{(1,2)}{(1,1)}{(2,2)}{(2,1)}}
,{\htwist~{(0,1)}{(0,0)}{(1,1)}{(1,0)}}
,{\vtwist~{(0,0)}{(0,-1)}{(0,0)}{(0,-1)}}
,{\vtwist~{(1,0)}{(1,-1)}{(1,0)}{(1,-1)}}
,{\vtwist~{(2,1)}{(2,-1)}{(2,1)}{(2,-1)}}
,{\htwist~{(3,3)}{(3,-1)}{(3,3)}{(3,-1)}}
,{\htwist~{(4,3)}{(4,-1)}{(4,3)}{(4,-1)}}
,{\htwist~{(5,3)}{(5,-1)}{(5,3)}{(5,-1)}}
\endxy
\end{picture}
\end{minipage}
}
\subfigure{%
\label{fig:braid2}%
\begin{minipage}[t]{0.35\textwidth}
\begin{picture}(4,2.5)(-3.7,0.2)
\xy
0;/r1pc/:
,{\vtwist~{(0,5)}{(0,3)}{(0,5)}{(0,3)}}
,{\vtwist~{(1,5)}{(1,4)}{(1,5)}{(1,4)}}
,{\htwist~{(3,5)}{(3,4)}{(2,5)}{(2,4)}}
,{\vtwist~{(4,5)}{(4,4)}{(4,5)}{(4,4)}}
,{\vtwist~{(5,5)}{(5,3)}{(5,5)}{(5,3)}}
,{\htwist~{(2,4)}{(2,3)}{(1,4)}{(1,3)}}
,{\htwist~{(4,4)}{(4,3)}{(3,4)}{(3,3)}}
,{\htwist~{(0,3)}{(0,2)}{(1,3)}{(1,2)}}
,{\htwist~{(2,3)}{(2,2)}{(3,3)}{(3,2)}}
,{\htwist~{(4,3)}{(4,2)}{(5,3)}{(5,2)}}
,{\vtwist~{(0,2)}{(0,0)}{(0,2)}{(0,0)}}
,{\htwist~{(1,2)}{(1,1)}{(2,2)}{(2,1)}}
,{\htwist~{(3,2)}{(3,1)}{(4,2)}{(4,1)}}
,{\htwist~{(2,1)}{(2,0)}{(3,1)}{(3,0)}}
,{\htwist~{(1,1)}{(1,0)}{(1,1)}{(1,0)}}
,{\htwist~{(4,1)}{(4,0)}{(4,1)}{(4,0)}}
,{\htwist~{(5,2)}{(5,0)}{(5,2)}{(5,0)}}
\endxy
\end{picture}
\end{minipage}
}
\caption{{The braids $\rho_0$ and $\rho_1$, for $p=3$}}
\label{fig:braids}
\end{figure}

Define braids $\rho_0,\rho_1 \in B_{2p}$ by
\[
\rho_0 = \sigma_{p-1}\sigma_{p-2} \cdots \sigma_1
\quad \text{and} \quad
\rho_1 =
\tau^{-1} \sigma_1 \sigma_3 \cdots \sigma_{2p-3} \sigma_{2p-1} \tau,
\]
where
\begin{equation} \label{eq:tau}
\tau=(\sigma_2 \sigma_4 \cdots \sigma_{2p-2})
(\sigma_3 \sigma_5 \cdots \sigma_{2p-3})
\cdots (\sigma_{p-2}\sigma_p\sigma_{p+2})
(\sigma_{p-1}\sigma_{p+1})(\sigma_p),
\end{equation}
see Figure \ref{fig:braids}.  The braids $\rho_i$ are obtained from
the ``monomial braids'' of \cite{mono} by deleting the central strand,
corresponding to the hyperplane $H_3=\ker(x_3)$ in the full monomial
arrangement, but not in the monomial deletion.  As in \cite{mono}, the
braid monodromy $\eta\co  G(\B) \to P_{{2p}}$ may be expressed
in terms of these braids, as follows.

Define pure braids
$Z_1^{},Z_2^{},A_{1,2}^{(1)},\dots,A_{1,2}^{(p)}$ in $P_{2p}$ by
$Z_1^{}=\rho_0^p$, $Z_2^{}=\rho_1^{} \rho_0^p \rho_1^{-1}$, and
$A_{1,2}^{(r)}=\rho_0^{r-p}\rho_1^2 \rho_0^{p-r}$ for $1 \le r \le p$.
Let $\gamma_j$ and $\gamma_{12:r}$ be meridian loops in $M(\B)$ about
the lines $H_j=\ker(x_j)$ and $H_{12:r}=\ker(x_1-\xi^r x_2)$.  These
loops generate the fundamental group $G(\B)$.

\begin{prop}
\label{prop:BraidMono}
The braid monodromy $\eta:G(\B) \to P_{2p}$ of the fiber bundle
$M(\A_p) \to M(\B)$ is given by $\eta(\gamma_j)=Z_j^{}$,
$\eta(\gamma_{12:r}) = A_{1,2}^{(r)}$.
\end{prop}

\begin{cor}
\label{cor:BraidMonoPres}
The fundamental group of $M(\A_p)$ has presentation
\[
G(\A_p) = \left\langle
\begin{aligned}
&\gamma_1, \gamma_2, \gamma_{12:1}, \dots, \gamma_{12:p}\\
&\quad y_1, y_2, y_3, \dots, y_{2p}
\end{aligned}
{\ } \bigg{|} {\ }
\begin{aligned}
\gamma_j^{-1} y_i^{} \gamma_j^{} &= \eta(\gamma_j^{})(y_i^{}) \\
\gamma_{12:r}^{-1} y_i^{} \gamma_{12:r}^{}
&= \eta(\gamma_{12:r}^{})(y_i^{}) \\
\end{aligned}
\right\rangle,
\]
where $i=1,\dots,2p$, $j=1,2$, and $r=1,\dots,p$, and the pure
braids $\eta(\gamma)$ act on the free group $\FF_{2p}=\langle
y_1,\dots,y_{2p} \rangle$ by the Artin representation.
\end{cor}

\subsection{Fundamental group of the decone}
\label{subsect:pi1dec}
Let $\Gamma=\gamma_1\gamma_{12:1}\cdots \gamma_{12:p-1}\gamma_2
\gamma_{12:p} \in G(\B)$.  Note that $\eta(\Gamma)=A_{[2p]} $ is
the full twist on all strands.  As is well known, this braid
generates the center of $P_{2p}$.  It follows that
$\Gamma$ is central in $G(\B)$, so
\[
G(\B)= \FF_{p+1} \times \Z=\langle \gamma_1,
\gamma_{12:1},\dots,\gamma_{12:p}\rangle
\times \langle \Gamma \rangle.
\]

To simplify calculations in \S\ref{sec:MoreHomology} below, we will
work with an explicit decone of the arrangement $\A_p$, as opposed to
the projectization.  Let $\dA_p$ denote the decone of $\A_p$ with
respect to the hyperplane $H_2=\ker(x_2)$.  This is an affine
arrangement in $\C^2$ (with coordinates $x_1,x_3$), defined by
$Q(\dA_p) = x_1(x_1^p-1)(x_1^p-x_3^p)(1-x_3^p)$.  From the above
discussion, we obtain the following presentation for the fundamental
group of the complement of $\dA_p$:
\begin{equation}
\label{eq:gdap}
G(\dA_p) = \left\langle
\begin{aligned}
&\gamma_1, \gamma_{12:1}, \dots, \gamma_{12:p}\\
&\quad y_1, y_2, \dots, y_{2p}
\end{aligned}
{\ } \bigg{|} {\ }
\begin{aligned}
\gamma_1^{-1} y_i^{} \gamma_1^{} &= \eta(\gamma_1^{})(y_i^{})\\
\gamma_{12:r}^{-1} y_i^{} \gamma_{12:r}^{}
&= \eta(\gamma_{12:r}^{})(y_i^{})
\end{aligned}
\right\rangle,
\end{equation}
where, as before,  $i=1,\dots,2p$ and $r=1,\dots,p$.

\subsection{Characteristic varieties}
\label{subsec:cvdel}
Set $n=3p+2=\abs{\A_p}$.  Denote the coordinates
of the algebraic torus $(\K^{\cross})^n$ by
$
z_1,z_2,z_{12:1},\dots,z_{12:p}, z_{13:1},
\dots,z_{13:p},z_{23:1},\dots,z_{23:p},
$
where $z_i$ corresponds to the hyperplane $H_i=\ker(x_i)$ and
$z_{ij:r}$ to the hyperplane $H_{ij:r}=\ker(x_i - \xi^r x_j)$.

The following theorem was proved for $p=2$ in \cite{suciu02},
and for $p\ge 2$ in \cite{triples}, in the case $\K=\C$.  The same
proofs work for an arbitrary, algebraically closed field $\K$.

\begin{thm}
\label{thm:char0CJL}
In addition to components of dimension $2$ or higher,
the variety $\Sigma^1_1(M(\A_p),\K)$ has $1$-dimensional components
$C_1,\dots, C_{p-1}$, given by
\[
\bigcup_{i=1}^{p-1} C_{i} = \Bigg\{
\left(u^p,v^p,w,\dots,w,v,\dots,v,u,\dots,u\right) \in (\K^{\cross})^n \:
\bigg| \bigg. \: \begin{aligned} 
&\tsum_{j=0}^{p-1} w^j = 0 \\ 
&\text{and}\ \: uvw=1
\end{aligned}
\Bigg\},
\]
where $C_i$ is obtained by setting $w$ equal to the $i$-th power of a 
fixed primitive $p$-th root of unity in $\K$.
\end{thm}

If $\ch\K=p$, then $C_i$ is a subtorus of $(\K^\cross)^n$, so passes
through the origin $\b{1}$.  However, if $\ch\K\neq p$, then $C_i$ is
a subtorus translated by a character of order $p$.  The results of
\S\ref{subsec:tj} imply that there exist choices of multiplicities
$\b{a}$ for $\A_{p}$ such that the first homology group of the
corresponding Milnor fiber, $F_{\b{a}}$, has $p$-torsion.  In
particular, we have the following.

\begin{cor}
\label{cor:FpTorsion}
Let $F_p=f_p^{-1}(1)$ be the Milnor fiber of the polynomial defined in
\eqref{eq:fp}.  Then $H_1(F_p;\Z)$ has $p$-torsion.
\end{cor}
\begin{proof}
Let $U=M(\dA_p)$ be the complement of the projectivization of
$\A_p$.  Note
that for $\b{t} \in C_i$, we have $\tau(\b{t})=1$.  So $C_i \subset
\Sigma_1^1(U,\K)$ by Proposition~\ref{prop:cvMandU}.  In the case
$\K=\C$, let $\b{s}_i=
\bigl(1,1,\xi^i,\dots,\xi^i,\xi^{-i},\dots,\xi^{-i},1,\dots,1\bigr)$.
where $\xi=\exp(2\pi\ii/p)$, and
\[
T=\{\bigl(u^p,v^{p},1,\dots,1,v,\dots,v,u,\dots,u \bigr)
\in(\C^\cross)^n
\mid uv=1\}.
\]
Then $\ord(\b{s}_i)=p$, $T$ is a one-dimensional subtorus of
$(\C^\cross)^n$, and $C_i= \b{s}_i T$.  One can check that $T
\not\subseteq \Sigma_1^1(U,\C)$ using known properites of
characteristic varieties of arrangements, see \cite{LY}.  Hence,
Theorem \ref{th:transl} implies that there are choices of
multiplicities $\b{a}$ for $\A_p$ for which $H_1(F_\b{a};\Z)$ has
$p$-torsion.  Arguing as in the proof of that theorem, and using Lemma
\ref{lem:mults}, reveals that among these choices of multiplicities
are $\b{a}=(2,1,3,3,2,2,1,1)$ in the case $p=2$, and
$\b{a}=(1,1,2,\dots,2,1,\dots,1,1,\dots,1)$ in the case $p\neq 2$.
These choices yield the polynomials $f_p$ of \eqref{eq:fp}.
\end{proof}

\section{Homology calculations}
\label{sec:MoreHomology}
Keeping the notation from the previous section, we analyze the
homology of $G(\dA_p)=\pi_1(M(\dA_p))$ with coefficients in the rank
one local systems that arise in the study of the Milnor fibration
$f_p\co  M(\A_p) \to\C^*$.  In this section, we consider the
case where $p \neq 2$ is an odd prime.

Let $\K$ be an algebraically closed field.  Recall that $\dA_p$ is the
decone of $\A_p$ with respect to the hyperplane $H_2=\ker(x_2)$, which
has multiplicity $1$ in the multi-arrangement defined by $f_p$.
Consequently, to analyze the homology of the Milnor fiber $F_p$ using
Theorem \ref{thm:hommf}, we will consider the modules $\K_{\b{t}(k)}$
corresponding to characters $\b{t}(k)$ defined by
\begin{equation}
\label{eq:tk}
\b{t}(k)=(t,t^2,\dots,t^2,t,\dots,t,t,\dots,t)
\in(\K^{\cross})^{n-1},
\end{equation}
where $t=\zeta^{N/k}$ is a power of a primitive $N$-th root of unity,
$N=4p+2$, $k\neq 1$ is a positive integer dividing $N$, and $n=3p+2$.

\begin{prop}
\label{prop:not2k}
If $k \ne 2$ and $\ch\K\nmid N$, then $H_1(G(\dA_p);\K_{\b{t}(k)})=0$.
\end{prop}
\begin{proof}
The braid $Z_1=\eta(\gamma_1)$ is a full twist on strands $1$ through
$p$, given in terms of the standard generators $A_{i,j}$ of $P_{2p}$ by
\[
Z_1=A_{1,2}(A_{1,3}A_{2,3}) \cdots \cdots (A_{1,p}\cdots
A_{p-1,p}).
\]
Consider the generating set $\{u_1,\dots,u_p,v_1,\dots,v_p\}$ for the
free group $\FF_{2p}$ given by $u_r = y_1y_2\cdots y_r$ and
$v_r=y_{p+r}$, $1\le r \le p$.  The action of the braid $Z_1$ on this
generating set is given by $Z_1^{}(u_i^{}) = u_p^{} u_i^{} u_p^{-1}$
and $Z_1^{}(v_j^{})=v_j^{}$ for $1\le i,j \le p$, see
\cite[\S6.4]{CScc}.

Taking $\gamma_1,\gamma_{12:r},u_r,v_r$ ($1\le r\le p$) as generators
for $G(\dA_p)$, we obtain from \eqref{eq:gdap} a presentation with
relations
\begin{align*}
u_i \gamma_1u_p &=\gamma_1 u_p u_i, & u_p \gamma_1&=\gamma_1 u_p,  &
v_j\gamma_1 &= \gamma_1 v_j, \\
u_j \gamma_{12:r}&=\gamma_{12:r}A_{1,2}^{(r)}(u_j), &
v_j \gamma_{12:r}&=\gamma_{12:r}A_{1,2}^{(r)}(v_j),
\end{align*}
where $1 \le i \le p-1$, $1\le j \le p$, and $1\le r \le p$.

Let $\b{A}$ denote the Alexander matrix obtained from this
presentation by taking Fox derivatives and abelianizing.  This is a
$2p(p+1) \times (3p+1)$ matrix with entries in the ring of Laurent
polynomials in the variables $\gamma_1,\gamma_{12:r},u_r,v_r$, and has
the form
\begin{equation*}
\label{eq:AlexanderMatrix}
\b{A} =
\begin{pmatrix}
\Delta & 0 & \cdots & 0 & 0 &
\II_{2p} - \gamma_1 \Theta(Z_1)\hfill  \\
0 & \Delta &  & 0 & 0 &
\II_{2p} - \gamma_{12:1} \Theta(A_{1,2}^{(1)})\hfill \\
\vdots & & \ddots & & & \phantom{\II_{2p}}\ \vdots\hfill \\
0 & 0 & & \Delta & 0 &
\II_{2p} - \gamma_{12:p-1} \Theta(A_{1,2}^{(p-1)}) \\
0 & 0 & \cdots & 0 & \Delta &
\II_{2p} - \gamma_{12:p} \Theta(A_{1,2}^{(p)})\hfill
\end{pmatrix},
\end{equation*}
where $\Delta$ is the column vector
$\left( u_1-1 , \dots , u_p-1 ,v_1-1,\dots,  v_p-1 \right)^{\top}$,
$\II_m$ is the $m\times m$
identity matrix, and $\Theta\co  P_{2p}\to \GL(2p,
\Z[ u_1^{\pm 1},\dots,u_p^{\pm 1},v_1^{\pm 1},\dots,v_p^{\pm 1}])$
is the Gassner representation.

Let $\b{A}(k)$ denote the evaluation of the Alexander matrix at the
character $\b{t}(k)$.  This evaluation is given by $\gamma_1\mapsto
t$, $\gamma_{12:r}\mapsto t^2$, $y_i\mapsto t$, so $u_r \mapsto t^r$
and $v_r\mapsto t$.  To show that $H_1(G(\dA_p);\K_{\b{t}(k)})=0$, it
suffices to show that $\b{A}(k)$ has rank $3p$.

A calculation (compare \cite[Prop.~6.6]{CScc}) reveals that the
evaluation at $\b{t}(k)$ of $\II_{2p} - \gamma_1 \Theta(Z_1)$ is upper
triangular, with diagonal entries $1-t^{p+1}$ and $1-t$.  Recall that
$\zeta$ is a primitive $N$-th root of unity, where $N=4p+2$, that $k
\neq 1$ divides $N$, and that $t=\zeta^{N/k}$.  Since $p$ is prime and
$k\neq 2$ by hypothesis, $k$ does not divide $p+1$.  Consequently, all
of the diagonal entries of the evaluation at $\b{t}(k)$ of $\II_{2p} -
\gamma_1 \Theta(Z_1)$ are nonzero.  It follows that $\rank
\b{A}(k)=3p$.
\end{proof}

For the character $\b{t}(2)=(-1,1,\dots,1,-1,\dots,-1,-1,\dots,-1)$,
and the corresponding module $\K_{\b{t}(2)}$, there are several cases
to consider.

First, note that if $\ch\K=2$, then $\b{t}(2)=\b{1}$ is the trivial
character.

If $\ch\K=p$, Theorem \ref{thm:char0CJL} and equation
\eqref{eq:decone} combine to show that $\b{t}(2) \in
\Sigma^1_1(M(\dA_p),\K)$.  Moreover, $\b{t}(2)\ne \b{1}$, since $p\neq
2$.  Hence, in this case the depth of $\b{t}(2)$ is at least $1$.

If $\ch\K\ne 2$ or $p$, one can show that the character $\b{t}(2)$
does not
lie in any component of $\Sigma^1_1(M(\dA_p),\K)$ of positive
dimension.  However, this does not rule out the possibility that
$\b{t}(2)$ is an isolated point in $\Sigma^1_1(M(\dA_p),\K)$.  This is
not the case, as the next result shows.

\begin{prop}
\label{prop:is2k}
Let $\K$ be an algebraically closed field.  If $\ch\K=p$, then
$\depth_\K(\b{t}(2))=1$.  If $\ch\K\neq 2$ or $p$, then
$\depth_\K(\b{t}(2))=0$.
\end{prop}

We will sketch a proof of this proposition by means of a sequence of
lemmas.  As above, we will analyze the Alexander matrix arising from a
well-chosen presentation of the group $G(\dA_p)$.

The presentation of $G(\dA_p)$ given in \eqref{eq:gdap} is
obtained from the realization of this group as a semidirect product,
$G(\dA_p) = \FF_{2p} \rtimes_{\overline{\alpha}} \FF_{p+1}$.  The
homomorphism $\overline{\alpha}\co  \FF_{p+1} \to \Aut(\F_{2p})$ is
the composition of the Artin representation with the braid monodromy
$\overline{\eta}\co  \FF_{p+1} \to P_{2p}$ given by
$\overline{\eta}\co  \gamma_1 \mapsto Z_1$, $\gamma_{12:r}\mapsto
A_{1,2}^{(r)}$.  We first modify the map $\overline{\eta}$, as
follows.

Recall the braid $\tau\in B_{2p}$ from \eqref{eq:tau}.  Conjugation by
$\tau$ induces an automorphism $\con_\tau\co  P_{2p} \to P_{2p}$,
$\beta \mapsto \tau\beta\tau^{-1}$.  Then, $\con_\tau \circ
\overline{\eta}\co  \FF_{p+1}\to P_{2p}$ is another choice of braid
monodromy for the (fiber-type) arrangement $\dA_p$, and the
presentation of $G(\dA_p)$ resulting from composing $\con_\tau\circ
\overline{\eta}$ and the Artin representation is equivalent to that
obtained from $\overline{\alpha}$.

\begin{lem}
\label{lem:combing}
In terms of the standard generating set for the pure braid group
$P_{2p}$, the braids $\con_\tau(A_{1,2}^{(p)})$ and $\con_\tau(Z_1)$
are given by
\begin{align*}
\con_\tau\big(A_{1,2}^{(p)}\big)&=A_{1,2}A_{3,4}\cdots A_{2p-1,2p},\\
\con_\tau(Z_1) &= A_{1,3}(A_{1,5}A_{3,5})\cdots
(A_{1,2p-1}A_{3,2p-1}\cdots A_{2p-3,2p-1}).
\end{align*}
\end{lem}
\begin{proof}
Recall that $\rho_0=\sigma_{p-1} \sigma_{p-2}\cdots \sigma_1$,
$\rho_1=\tau^{-1}\sigma_1\sigma_3\cdots\sigma_{2p-1}\tau$, and
$A_{1,2}^{(r)}=\rho_0^{r-p}\rho_1^2\rho_0^{p-r}$.  Hence,
$\con_\tau(A_{1,2}^{(p)})=\sigma_1^2\sigma_3^2\cdots\sigma_{2p-1}^2=
A_{1,2}A_{3,4}\cdots A_{2p-1,2p}$.  Also recall that $Z_1=A_{[p]}=
A_{1,2}(A_{1,3}A_{2,3})\cdots \cdots(A_{1,p}A_{2,p}\cdots A_{p-1,p})$
is the full twist on strands $1$ through $p$.  We will show that
$\con_\tau(Z_1)$ is as asserted (for any integer $p\ge 2$) by
induction on $p$.

Write $\tau=\tau_p$.  When $p=2$, we have
$\tau_2^{}A_{[2]}^{}\tau^{-1}_2=\sigma_2^{}
\sigma_1^2\sigma_2^{-1}=A_{1,3}$.  So inductively assume that
$\tau_p^{}A_{[p]}^{}\tau_p^{-1}=A_{O[p]}^{}$, where
$O[p]=\{1,3,\dots,2p-1\}$.  Using \eqref{eq:tau} and the braid
relations, we have $\tau_{p+1}=\tau_p \beta_{p+1}$, where
$\beta_{p+1}=\sigma_{2p}\sigma_{2p-1} \cdots \sigma_{p+1}$.  Note that
$\beta_{p+1}$ commutes with $A_{[p]}$.  Hence,
\[
\begin{aligned}
\tau_{p+1}^{}A_{[p+1]}^{}\tau_{p+1}^{-1}
&=\tau_p\beta_{p+1}A_{[p]}(A_{1,p+1}A_{2,p+1}\cdots
A_{p,p+1})\beta_{p+1}^{-1}\tau_p^{-1}\\
&=\tau_pA_{[p]}\tau_p^{-1}\cdot\tau_p\beta_{p+1}(A_{1,p+1}A_{2,p+1}
\cdots A_{p,p+1})\beta_{p+1}^{-1}\tau_p^{-1}\\
&=A_{O[p]} \cdot \tau_p\beta_{p+1}(A_{1,p+1}A_{2,p+1}\cdots
A_{p,p+1})\beta_{p+1}^{-1}\tau_p^{-1}\\
&=A_{O[p]} \cdot \tau_p(A_{1,2p+1}A_{2,2p+1}\cdots
A_{p,2p+1})\tau_p^{-1}
\end{aligned}
\]
by induction, and the readily checked fact that
$\beta_{p+1}A_{i,p+1}\beta_{p+1}^{-1}=A_{i,2p+1}$.
The result now follows from the equality
$\tau_p A_{i,2p+1}\tau_p^{-1}=A_{2i-1,2p+1}$, which
may itself be established by an inductive argument.
\end{proof}

Write $\mathfrak{z}=\con_\tau(Z_1)$ and
$\mathfrak{a}=\con_\tau(A_{1,2}^{(p)})$.  We specify a generating set
for the free group $\FF_{2p}=\langle y_1,\dots,y_{2p}\rangle$ for
which the action of these braids is tractable.  For $1\le r \le p$,
let $u_r=y_1y_2 \cdots y_{2r-1}y_{2r}$ and $v_r=y_{2r-1}$.  Write
$V=v_1v_2\cdots v_p$.  It is readily checked that the set of elements
$\{u_1,\dots,u_p,v_1,\dots,v_p\}$ generates $\FF_{2p}$.  Moreover, a
calculation using the Artin representation yields the following.

\begin{lem}
\label{lem:zact}
The action of the braids $\mathfrak{z}$ and
$\mathfrak{a}$ on the set $\{u_r,v_r\}_{r=1}^p$ is given by
\[
\begin{matrix}
\mathfrak{z}(u_r) = u_r [v_{r+1}\cdots v_p,v_1\cdots v_r], &
\mathfrak{a}(u_r)=u_r, \hfill \\
\mathfrak{z}(v_r) =Vv_rV^{-1}, \hfill &
\mathfrak{a}(v_r)=u_{r-1}^{-1}u_r^{}v_r^{} u_r^{-1}u_{r-1}^{}.\\
\end{matrix}
\]
\end{lem}

Note that $\mathfrak{z}(u_p)=u_p$ and that
$\mathfrak{a}(v_1)=u_1v_1u_1^{-1}$.

Now consider the presentation of the group $G(\dA_p)$ obtained from
the braid monodromy $\con_\tau \circ \, \overline{\eta}\co  \FF_{p+1}
\to P_{2p}$ and the Artin representation, using the generating set
$\{u_r,v_r\}_{r=1}^p$ for the free group $\FF_{2p}$.  Identify the
generators $\gamma_1,\gamma_{12:j}$ of $\FF_{p+1}$ with their images
in $P_{2p}$ via $\gamma_1 \mapsto \con_\tau(\gamma_1)=\mathfrak{z}$,
$\gamma_{12:p}\mapsto \con_\tau(\gamma_{12:p})=\mathfrak{a}$, and
write $\gamma_{12:j}\mapsto \con_\tau(\gamma_{12:j})=\mathfrak{a}_j$
for $1\le j \le p-1$.  With this notation, the presentation for
$G(\dA_p)$ has relations
\begin{equation}\label{eq:UsefulRelations}
\begin{matrix}
u_r\mathfrak{z} = \mathfrak{z} u_r [v_{r+1}\cdots v_p,v_1\cdots v_r],
& u_r\mathfrak{a}=\mathfrak{a}u_r, \hfill \\
v_r\mathfrak{z}=\mathfrak{z}Vv_rV^{-1}, \hfill & v_r\mathfrak{a}=
\mathfrak{a}u_{r-1}^{-1}u_r^{}v_r^{} u_r^{-1}u_{r-1}^{}, \hfil\\
\end{matrix}
\end{equation}
and $u_r^{}\mathfrak{a}_j^{}=
\mathfrak{a}_j^{}\cdot\mathfrak{a}_j^{}(u_r^{})$,
$v_r^{}\mathfrak{a}_j^{}=
\mathfrak{a}_j^{}\cdot\mathfrak{a}_j^{}(v_r^{})$,
for $1\le j \le p-1$ and $1\le r \le p$.

Let $\b{A}$ be the Alexander matrix obtained from this presentation,
and $\b{A}(2)$ the evaluation at the character $\b{t}(2)$.  This
evaluation is given by $\gamma_1\mapsto -1$, $\gamma_{12:j} \mapsto
1$, $y_r \mapsto -1$, so $\mathfrak{z}\mapsto -1$, $\mathfrak{a}
\mapsto 1$, $\mathfrak{a}_j \mapsto 1$, $u_r \mapsto 1$,
$v_r \mapsto -1$.  Let $\sM$, $\sA$, and $\sA_j$ denote the evaluations
at $\b{t}(2)$ of the Fox Jacobians of the actions of the pure braids
$\mathfrak{z}$, $\mathfrak{a}$, and $\mathfrak{a}_j$, respectively.
With this notation, we have
\begin{equation*}
\label{eq:A2Matrix}
\b{A}(2)=
\text{\small $
\begin{pmatrix}
\Delta(2) & 0 & \cdots & 0 & 0 & \II_{2p} +\sM\hfill  \\
0 & \Delta(2) &   & 0 & 0 &  \II_{2p} - \sA_1\hfill \\
\vdots & & \ddots & & &  \phantom{\II_{2p}}\ \vdots\hfill \\
0 & 0 &  & \Delta(2) & 0 & \II_{2p} - \sA_{p-1} \\
0 & 0 & \cdots & 0 & \Delta(2) &  \II_{2p} - \sA\hfill
\end{pmatrix}$},
\end{equation*}
where $\Delta(2)=\begin{pmatrix}0 & \cdots & 0 & -2 & \cdots & -2
\end{pmatrix}^\top$ is the evaluation of $\Delta$ at $\b{t}(2)$.
Note that the entries of $\b{A}(2)$ are integers, and recall that
$\b{A}(2)$ has size $2p(p+1) \times (3p+1)$.

To establish Proposition \ref{prop:is2k}, we must show that
$\rank_\K \b{A}(2)=3p-1$ or $3p$, according to whether the field $\K$
has characteristic $p$ or not (recall that  $\ch\K\neq 2$ by 
assumption).
In the case $\ch\K=p$, we already know that
$\b{t}(2)$ belongs to $\Sigma^1_1(M(\dA_p),\K)$, so the inequality
$\rank_\K \b{A}(2) \le 3p-1$ holds.  Thus, it suffices to prove the
next result.

\begin{lem}
\label{lem:snf}
The (integral) Smith normal form of the matrix $\b{A}(2)$
has diagonal entries $2,\dots,2$ (repeated $3p-1$ times) and $2p$.
\end{lem}
\begin{proof}
The matrix $\b{A}(2)$ is equivalent, via row and
column operations, to the matrix
\begin{equation}
\label{eq:forSNF}
\text{\small $
\begin{pmatrix}
0 & 0 & \cdots & 0 & 0 &  2\II_{2p} +\sM- \sA  \\
0 & \Delta(2) &  & 0 & 0 &  \phantom{2}\II_{2p} - \sA_1\hfill \\
\vdots & & \ddots & & &\phantom{2\II_{2p}}\ \vdots\hfill \\
0 & 0 &  & \Delta(2) & 0 &  \phantom{2}\II_{2p} -
\sA_{p-1}\hfill \\
0 & 0 & \cdots & 0 & \Delta(2)&  \phantom{2}\II_{2p} - \sA\hfill
\end{pmatrix}$}.
\end{equation}
A Fox calculus exercise using \eqref{eq:UsefulRelations} shows that
all entries of the matrices $\II_{2p}-\sA$ and $\II_{2p}-\sA_j$,
$1\le j\le p-1$, are divisible by $2$, and that
\[
2\II_{2p} +\sM- \sA =
2\begin{pmatrix}
\II_p & {\sf P}\\ {\sf L}-\II_p & {\sf Q}
\end{pmatrix},
\]
where ${\sf L}_{i,j}=\delta_{i,j+1}$ (Kronecker delta),
${\sf Q}_{i,j}=(-1)^{j+1}$, and
\[
{\sf P}_{i,j}=\begin{cases}
(-1)^{j} & \text{if $i$ odd and $j>i$},\\
(-1)^{j+1} & \text{if $i$ even and $j\le i$},\\
0 & \text{otherwise}.
\end{cases}
\]
Let ${\sf U}=\begin{pmatrix}
\II_p & 0 \\ \II_p - {\sf L} & \II_p\end{pmatrix}
\begin{pmatrix} \II_p & 0 \\ 0 & {\sf R}\end{pmatrix}$ and
${\sf V}=\begin{pmatrix} \II_p & -{\sf P}\\ 0 & \pmi \II_p
\end{pmatrix}
\begin{pmatrix} \II_p & 0 \\ 0 & {\sf S}\end{pmatrix}$, where
\[
\sf R=
\text{\small $
\begin{pmatrix}
1 & 0 & 0 &  \cdots & 0 & \pmi 0 \\
0 & 1 & 0 &  & 0 &-1\\
0 & 2 & 1 &  & 0 &-2\\
\vdots & & & \ddots  &  & \pmi \vdots \\
0 & 2 & 2 &   & 1 & -(p-2)\\
0 & 2 & 2 &  \cdots & 2 & -(p-2)
\end{pmatrix}$},\quad
\sf S=
\text{\small $
\begin{pmatrix}
1 & 0 & 0 &  \cdots & 0 & \pmi 0 \\
0 & 1 & 0 &   & 0 &-1\\
0 & 0 & 1 &   & 0 &-2\\
\vdots &  & &\ddots  &  & \pmi \vdots \\
0 & 0 & 0 & & 1 & -(p-2)\\
0 & 0 & 0 & \cdots & 0 & \pmi 1
\end{pmatrix}$}.
\]
Then one can check that $\det{\sf R}=\det{\sf S}=1$, and that
${\sf U}(2\II_{2p} +\sM- \sA){\sf V}$ is a $2p\times 2p$ diagonal
matrix with diagonal entries $2,\dots,2,2p$ (in this order).
Using these facts, further row and
column operations reduce the matrix \eqref{eq:forSNF} to
\begin{equation*}
\begin{pmatrix}
0 & 2 \II_{3p-1} & 0  \\
0 & 0            & 2p \\
0 & 0            & {\sf v}
\end{pmatrix},
\end{equation*}
where ${\sf v}$ is a column vector whose entries are even integers.
Now recall that if $\K$ is a field of characteristic $p$, then
$\rank_\K \b{A}(2) \le 3p-1$.  Consequently, the entries of $\sf v$
must be divisible by $p$.  The result follows.
\end{proof}

\section{Proof of Theorem~\ref{theo:h1fp}}
\label{sec:proofthm}

We are now in position to complete the proof of
Theorem~\ref{theo:h1fp} from the Introduction.  Recall we are given a
prime $p$ and the homogenous polynomial $f_p$ specified in
\eqref{eq:fp}, and we need to compute the first homology group of the
Milnor fiber $F_p=f_p^{-1}(1)$.  We shall treat the cases of odd and
even primes $p$ separately.

\subsection{The case $p\ne 2$}
\label{subsect:podd}
Recall that $\A_p$ is the arrangement in $\C^{3}$ defined by the
polynomial 
$Q(\A_p)=x_1^{}x_2^{}(x_1^p-x_2^p)(x_1^p-x_3^p)(x_2^p-x_3^p)$. 
The choice of multiplicities
\[
\b{a}=(1,1,2,\dots,2,1,\dots,1,1,\dots,1)
\]
yields the homogeneous polynomial
$f_p=x_1^{}x_2^{}(x_1^p-x_2^p)^2(x_1^p-x_3^p)(x_2^p-x_3^p)$.  This
gives rise to a Milnor fibration $f_p\co  M(\A_p)\to \C^*$, with
fiber $F_p=F_{\b{a}}$.  Let $\Z_N$ be the cyclic group of order
$N=\deg(f_p)=4p+2$, with generator $g$.  The $N$-fold cyclic cover
$F_p \to M(\dA_p)$ is classified by the epimorphism $\lambda\co 
G(\dA_p) \to \Z_N$ given by $\lambda(\gamma_1)=g$,
$\lambda(\gamma_{12:r})=g^2$, $\lambda(\gamma_{i3:r})=g$.

Let $\K$ be an algebraically closed field, of characteristic not
dividing $N$.  The homology group $H_1(F_p;\K)$ may be calculated
using Theorem~\ref{thm:hommf}:
\begin{equation*}
\label{eq:dimkfp}
\dim_\K H_1(F_p;\K) = 3p+1 + \sum_{1\neq k|N} \varphi(k)
\depth_\K \b{t}(k),
\end{equation*}
where $\b{t}(k)$  are the characters defined in \eqref{eq:tk}.
Using  Propositions \ref{prop:not2k} and \ref{prop:is2k}, we find:
\begin{equation}
\label{eq:dimh1}
\dim_\K H_1(F_p;\K) =
\begin{cases}
3p+1, &\text{if $\ch\K\nmid 2p(2p+1)$,}\\
3p+2 , &\text{if $\ch\K=p$.}
\end{cases}
\end{equation}

Now recall that we have an isomorphism $H_1(F_p;\Z)\cong
H_1(G;\Z[\Z_N])$ between the first homology of $F_{p}$ and that of
$G=G(\dA_{p})$, with coefficients in the $G$-module $\Z[\Z_{N}]$
determined by the epimorphism $\lambda\co  G \to \Z_{N}$.
Let $\Z_2\subset \Z_N$ be the subgroup generated by $g^{N/2}$,
and let $\Z[\Z_2]\subset \Z[\Z_N]$ be the corresponding $G$-submodule.
Denote by $J$ the kernel of the augmentation
map $\epsilon\co  \Z[\Z_2]\to \Z$. Notice that $g^{N/2}$ acts on
$J\cong \Z$ by multiplication by $-1$. Hence, the induced $G$-module
structure on $J$ is given by the composite
$G \xrightarrow{\alpha \circ \ab} \Z^{3p+1}
\xrightarrow{\b{t}(2)} \{\pm 1\}$,
which shows that $J$ is the integral analogue of the local system
$\K_{\b{t}(2)}$.  Let $Q=\Z[\Z_N]/J$ be the quotient $G$-module,
and consider the homology long exact sequence corresponding
to the coefficient sequence $0 \to J \to \Z[\Z_{N}] \to Q \to 0$:
\begin{equation}
\label{eq:les}
\cdots \to H_2(G;Q) \to H_1(G;J) \to H_1(G;\Z[\Z_N]) \to H_1(G;Q)
\to \cdots
\end{equation}

By Lemma~\ref{lem:snf}, we have
$H_1(G;J) \cong (\Z_2)^{3p}\oplus \Z_p$.
Over an algebraically closed field $\K$ with $\ch\K\nmid N$, the
$G$-module $Q$ decomposes as the direct sum of the modules
$\K_{\b{t}(k)}$, $k \neq 2$, together with the trivial module.  So
Proposition \ref{prop:not2k} implies that that $H_1(G;Q)$ has no
$q$-torsion, for any odd prime $q$ not dividing $2p+1$.  Note that
$H_2(G;Q)$ is free abelian, since the cohomological dimension of
$G=\b{F}_{2p} \rtimes_{\bar\alpha} \b{F}_{p+1}$ is $2$.  Applying
these observations to the long exact sequence \eqref{eq:les} reveals
that the map $H_1(G;J)\to H_1(G;\Z[\Z_N])$ induces an
isomorphism on $p$-torsion.  Therefore:
\begin{equation}
\label{eq:h1fp}
H_1(F_p;\Z)=\Z^{3p+1}\oplus \Z_p\oplus T,
\end{equation}
where $T$ is a finite abelian group such that $T\otimes \Z_q=0$ if
$q\nmid 2(2p+1)$.  This finishes the proof of Theorem~\ref{theo:h1fp}
in the case $p\ne 2$.

\begin{rem}
\label{lem:monop}
The $p$-torsion in \eqref{eq:h1fp} appears in the $(-1)$-eigenspace of
the algebraic monodromy $h_*$, see Lemma~\ref{lem:Fsplits}.  Since an
automorphism of $H_1(F_p;\Z)$ must preserve the $p$-torsion elements,
$h_*$ acts on the $\Z_{p}$ direct summand by $x\mapsto -x$.
\end{rem}

\begin{figure}[ht]
\setlength{\unitlength}{0.7cm}
\begin{picture}(4,4.2)(-6,-0.5)
\multiput(2.5,0)(0.5,0){3}{\line(0,1){3}}
\put(1.5,2){\line(1,0){3}}
\put(1.5,1){\line(1,0){3}}
\put(1.5,0){\line(1,1){3}}
\put(4.5,0){\line(-1,1){3}}
\put(4.8,-0.25){\makebox(0,0){$2$}}
\put(4.8,3.25){\makebox(0,0){$2$}}
\put(4.8,1){\makebox(0,0){$1$}}
\put(4.8,2){\makebox(0,0){$1$}}
\put(3.5,3.5){\makebox(0,0){$3$}}
\put(3,3.5){\makebox(0,0){$2$}}
\put(2.5,3.5){\makebox(0,0){$3$}}
\end{picture}
\caption{{Decone of deleted
$\operatorname{B}_3$ arrangement, with multiplicities}}
\label{fig:deletedb3}
\end{figure}
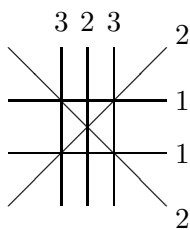

\subsection{The case $p=2$}
\label{subsect:peven}
Now consider the arrangement $\A_2$ in $\C^{3}$ defined by the
polynomial 
$Q(\A_2)=x_1^{}x_2^{}(x_1^2-x_2^2)(x_1^2-x_3^2)(x_2^2-x_3^2)$.  
This is a deletion of the $\operatorname{B}_3$ reflection arrangement,
and appears as Example~4.1 in \cite{suciu02} and Example~9.3 in
\cite{MShall}.  The polynomial $f_{2}=x_1^2x_2(x_1^2-x_2^2)^3
(x_1^2-x_3^2)^2(x_2-x_3)$ corresponds to the choice of multiplicities
$\b{a}=(2,1,3,3,2,2,1,1)$, shown in Figure~\ref{fig:deletedb3} (the
hyperplane at infinity has multiplicity $1$).

The variety $\Sigma^1_1(M(\A_2),\C)$ contains a $1$-dimensional
component $\b{s}T$, where $T=\set{(u^2,v^2,1,1,v,v,u,u) \mid uv=1}$
and $\b{s}=(1,1,-1,-1,-1, -1,1,1)$, see Theorem~\ref{thm:char0CJL}. 
The subtorus $T$ is not a component.  For example, the point $\b{t}\in
T$ given by $u=\exp(2\pi \ii/3)$ and $v=u^2$ is not in
$\Sigma^1_1(M(\A_2),\C)$.

The Milnor fiber $F_2=f_{2}^{-1}(1)$ is an $N$-fold cover of
$M(\dA_2)$, with $N=15$.  Using Theorem~\ref{thm:hommf} as before, we
find that $\dim_\K H_1(F_2;\K)=7$ if $\ch\K \ne 2,3$, or $5$, and
$\dim_\K H_1(F_2;\K)=9$ if $\ch\K =2$.  Direct computation with
the Alexander matrix of $G(\dA_2)$ (see \cite[Ex.~4.1]{suciu02})
gives the precise answer:
\begin{equation}
\label{eq:h1f2}
H_1(F_2;\Z)=\Z^7\oplus\Z_2\oplus\Z_2.
\end{equation}
This finishes the proof of Theorem \ref{theo:h1fp} in the remaining
case $p=2$.

\begin{rem}
\label{lem:mono2}
Once again, the monodromy action preserves the torsion part in
\eqref{eq:h1f2}, so $\Z_{15}$ acts on $\Z_2\oplus\Z_2$.  Since the
torsion in $H_1(F_2;\Z)$ appears in the eigenspaces of order $3$, the
monodromy acts via an automorphism of order $3$, which, in a suitable
basis, has matrix
$\bigl(\begin{smallmatrix}0&1\\1&1\end{smallmatrix}\bigr)$.
\end{rem}

\bibliographystyle{gtart}

\Addresses\recd

\end{document}